COMPRENDRE LES MATHÉMATIQUES POUR COMPRENDRE PLATON - *THÉÉTÈTE* (147d-148b)
*UNDERSTANDING MATHEMATICS TO UNDERSTAND PLATO* - THEAETEUS (147d-148b)

Salomon OFMAN
Institut mathématique de Jussieu-PRG
Histoire des Sciences mathématiques
4 Place Jussieu
75005 Paris
salomon.ofman@imj-prg.fr

**Résumé.** Dans cet article, nous étudions essentiellement les premières lignes de ce que l'on nomme traditionnellement la 'partie mathématique' du *Théétète* de Platon, où un jeune Athénien, Théétète, rapporte une leçon de mathématiques sur l'incommensurabilité de certaines grandeurs, à laquelle il a assisté. En termes modernes, il s'agit de la question de la rationalité (ou de l'irrationalité) des racines carrées des nombres entiers. En tant que le plus ancien texte qui nous soit parvenu sur le sujet, mais aussi sur les mathématiques et les mathématiciens grecs, sa valeur est inestimable. Les difficultés pour l'interpréter proviennent de l'étroite imbrication qu'on y trouve entre différents domaines : philosophie, histoire et mathématiques. Mais inversement, convenablement compris, il peut fournir des témoignages à la fois sur la question des origines de la théorie des irrationnels dans les mathématiques grecques et sur certains points de la pensée platonicienne. À partir d'une analyse mathématique prenant en compte le contexte historique et l'arrière-plan philosophique du dialogue généralement négligés, nous obtenons une interprétation nouvelle de ce texte qui, loin d'être un hommage à certains mathématiciens, est une critique radicale de leurs manières de penser. Et la leçon mathématique, loin d'être un hommage à de futurs succès mathématiques, apparaît, de manière cohérente avec le dialogue tout entier, se conclure sur une aporie.

**Abstract**. *In this paper, we study the so-called 'Mathematical part' of Plato's Theaetetus. Its subject concerns the incommensurability of certain magnitudes, in modern terms the question of the rationality or irrationality of the square roots of integers. As the most ancient text on the subject, and on Greek mathematics and mathematicians as well, its historical importance is enormous. The difficulty to understand it lies in the close intertwining of different fields we found in it: philosophy, history and mathematics. But conversely, correctly understood, it gives some evidences both about the question of the origins of the irrationals in Greek mathematics and some points concerning Plato's thought. Taking into account the historical context and the philosophical background generally forgotten in mathematical analyses, we get a new interpretation of this text, which far from being a tribute to some mathematicians, is a radical criticism of their ways of thinking. And the mathematical lesson, far from being a tribute to some future mathematical achievements, is ending on an* aporia*, in accordance with the whole dialogue.*

Mots-clés : aporie, définition, irrationnels (origine des), mathématiques, philosophie, Platon, science, *Théétète* (partie mathématique du).

## 1. Présentation de l'ouvrage : le prologue

Le dialogue proprement dit est précédé d'un prologue entre deux personnages Euclide[1] et Terpsion. Cela se passe à Mégare jusqu'où Euclide a raccompagné Théétète agonisant, après une bataille à proximité de la ville de Corinthe[2], non suite à des blessures reçues au combat, mais à une maladie, la dysenterie, qui a touché l'armé athénienne. Les deux personnages ayant beaucoup marché, décident de se reposer. Terpsion demande alors à son compagnon de lui exposer une discussion qui eut lieu peu avant la mort de Socrate. Ce dernier en avait fait le récit à Euclide qui s'était empressé de le mettre par écrit, le faisant corriger par Socrate 'toutes les fois qu'il avait l'occasion de le voir'[3]. Ce texte met en scène outre Socrate, Théétète, alors très jeune, et Théodore, un mathématicien de Cyrène, une colonie grecque située dans le nord-est de l'actuelle Lybie. Toutefois, dans le récit qui nous concerne ici, seuls les deux premiers participent directement au dialogue.

---

[1] Il s'agit d'Euclide de Mégare, longtemps confondu avec Euclide d'Alexandrie, le mathématicien auteur des *Éléments*.
[2] La datation de cette bataille est l'objet d'âpres disputes. Les deux dates communément citées sont soit autour de 390 BCE soit de 369 BCE ('*Before Common Era*').
[3] '… καὶ ὁσάκις Ἀθήναζε ἀφικοίμην, ἐπανηρώτων τὸν Σωκράτη ὃ μὴ ἐμεμνήμην' (142d7). On peut se demander si ce nombre ('ὁσάκις') est bien grand, cette rencontre objet du récit ayant eu lieu peu avant la condamnation à mort de Socrate (210d).

## 2. Introduction

En réponse à Socrate qui lui demande des nouvelles des jeunes Athéniens qui le fréquentent et suivent son enseignement mathématique, Théodore décrit de manière dithyrambique un de ses jeunes élèves, Théétète, qui est à la fois extrêmement doué et d'un naturel ouvert et généreux, mais laid, autant ajoute Théodore, que Socrate. À la demande de celui-ci, il appelle le jeune garçon qui sort précisément du gymnase. Ainsi débute une scène à trois interlocuteurs, Socrate, Théodore et Théétète.

Socrate, après avoir mis en doute les connaissances de Théodore sur la beauté et la laideur physiques, et amené, quoique avec réticence, Théétète à en convenir, propose un thème de débat : une discussion/définition ($\lambda\acute{o}\gamma o\varsigma$) sur ce qu'est la science/connaissance ($\grave{\epsilon}\pi\iota\sigma\tau\acute{\eta}\mu\eta$)[4].

Dans le rôle de répondant, Théodore se désistant, propose à sa place son élève aux si nombreuses qualités, Théétète.

Après beaucoup d'hésitations celui-ci accepte, et en réponse à la question de Socrate 'qu'est-ce que la science ?', en énumère une série, la géométrie, l'astronomie, l'harmonie (pour la musique), le calcul, déjà nommées par Socrate auxquelles il ajoute la cordonnerie et les autres techniques manuelles[5].

Socrate lui reproche sa 'générosité' (146d4-5)[6] qui, peu auparavant, avait été au contraire louée par Théodore (144d3), le jeune garçon donnant une pluralité d'exemples, là où on lui demandait une réponse, ce qu'est **la** science.

Suite à cette véhémente critique, Théétète se souvient soudainement d'un problème similaire, dit-il, que lui et un de ses camarades, homonyme de Socrate[7], s'étaient posé à la suite d'une leçon de Théodore à laquelle ils avaient assisté. C'est le début de ce qu'on appelle la 'partie mathématique du *Théétète*'.

---

[4] En fait, comme l'indique dès le début Socrate (144d8-145b1), cette enquête porte tout autant sur Théodore (et Théétète), son caractère, son témoignage, mais aussi, comme on le verra, sur ses mathématiques et son enseignement. La polysémie de termes fondamentaux dans le texte, en particulier '$\lambda\acute{o}\gamma o\varsigma$' et '$\grave{\epsilon}\pi\iota\sigma\tau\acute{\eta}\mu\eta$' rend délicate sa traduction mais également parfois sa compréhension.

[5] Mélange que Socrate rejettera immédiatement après, modifiant les 'techniques' de Théétète en 'sciences' (146d8).

[6] On a déjà ici un exemple d'incommensurabilité dont il sera question dans le passage mathématique. Dans l'introduction du *Politique*, où se retrouvent les mêmes personnages que ceux du *Sophiste*, dialogue censé se dérouler le lendemain de celui conté dans le *Théétète*, Socrate reprend vertement Théodore lorsque celui-ci met sur le même plan politique et philosophie, lui faisant convenir que ce sont des sujets qu'aucune sorte de proportion ne peut rendre (257a9-b8). Théétète met dans le même sac géométrie et cordonnerie, sciences ($\grave{\epsilon}\pi\iota\sigma\tau\hat{\eta}\mu\alpha\iota$) et techniques ($\tau\acute{\epsilon}\chi\nu\alpha\iota$) (146d7-8). Socrate le corrige sans le faire remarquer, comme il le fera à plusieurs reprises par la suite (en particulier en 148d3-6), en spécifiant la 'science de la cordonnerie' (146d7-8) et 'la science de la menuiserie' (148e1-2), remplaçant au passage le pluriel '$\tau\acute{\epsilon}\chi\nu\alpha\iota$' ('techniques') par le singulier '$\grave{\epsilon}\pi\iota\sigma\tau\acute{\eta}\mu\eta$' ('science') (146d7-8). La 'générosité' de Théétète le conduit non seulement à répondre à côté de la question et à remplacer l'unité par une multiplicité, mais également à mélanger des objets n'ayant aucune commune mesure. Diogène Laërce rapporte une anecdote analogue, mais dirigée cette fois contre la 'générosité' de Platon, c'est-à-dire son bavardage (Laërce 1965, VI, 2, 40).

[7] C'est ainsi qu'il est nommé par Platon. Suivant la tradition, il est appelé Socrate le Jeune pour de le distinguer du philosophe.

Afin d'en faciliter la compréhension, l'exposition mathématique que nous allons présenter sera parfois volontairement anachronique. Une telle transcription n'est pas exempte de dangers[8], car elle peut dissimuler des difficultés conceptuelles[9]. C'est pourquoi nous en ferons, pour chaque cas, une transcription dans les termes et les limites des mathématiques grecques de cette époque.

---

[8] Un exemple très simple est la commutativité de la multiplication (à savoir que pour deux entiers $m$ et $n$, le produit de $m$ par $n$ est égal à celui de $n$ par $m$), pure trivialité pour les modernes, alors que dans les *Éléments* d'Euclide, c'est la seizième proposition du livre VII, nécessitant le passage par la théorie des rapports des entiers, élaborée dans les propositions précédentes.

[9] Un exemple en est la démonstration d'irrationalité par 'anthyphérèse' (cf. infra, §5) qui est généralement donnée sur un unique un exemple et sous une forme algébrique et symbolique modernes (par exemple van der Waerden 1963, p. 143-146).

### 3. Les grandes lignes du récit

L'importance de la 'partie mathématique' du *Théétète* ne tient pas seulement en ce qu'elle est le plus ancien texte portant sur la mathématique grecque et les grandeurs irrationnelles, mais également sur le travail et la pratique des mathématiciens de l'Antiquité classique ($5^{\text{ème}}$ siècle BCE).

Le jeune Théétète commence son récit en rapportant une leçon donnée par Théodore à propos de l'incommensurabilité avec l'unité des côtés des carrés ayant pour aires *3* pieds, *5* pieds, jusqu'à *17* pieds[10]. En termes modernes, il s'agit de l'irrationalité de *√3*, *√5*, …, *√17*.
À la suite de ce cours, Théétète et son camarade Socrate se réunissent ensemble pour le retravailler. Ils affirment tout d'abord l'évidence d'une infinité de telles grandeurs. Aussi décident-ils d'en obtenir une caractérisation générale plutôt qu'au 'cas par cas' ('κατὰ μίαν ἑκάστην') à la manière de Théodore. Ils remarquent qu'on peut traduire cette question en termes géométriques en associant à chaque entier une figure plane, plus exactement des rectangles de surface ces entiers[11].
Ils subdivisent alors les rectangles, et donc les entiers, en deux classes. L'une est formée des carrés (de côtés entiers), ainsi ceux d'aire *9,…* , l'autre des rectangles qui ne peuvent être mis sous cette forme, ainsi ceux d'aire *5, 7, ….* . À tous les rectangles de cette seconde classe correspondent des grandeurs incommensurables avec celles de la première, à savoir les côtés des carrés ayant pour aire ces rectangles. En termes modernes, la première classe est formée des nombres entiers, la seconde d'irrationnels, les racines carrées des entiers non carrés parfaits[12].

À la différence de la leçon de Théodore, on a quitté la mesure physique des surfaces, donnée en pieds, l'énoncé des garçons concernant non plus de telles surfaces effectivement tracées, ou que l'on peut tracer, à la main, mais des nombres et des grandeurs abstraites. Le dessin, au moins partiellement, a été remplacé par une certaine forme de raisonnement. On peut donc estimer que cette mathématique est plus sophistiquée et plus proche de la mathématique telle que Platon la conçoit, par exemple au livre VII de la *République*, que celle dans la leçon de

---

[10] Théétète utilise ici une mesure de longueur ('ποδιαία') pour des surfaces. Certains traducteurs (ainsi Fowler 1996, McDowell 1973 ou plus récemment Chappell 2004) n'hésitent pas à le rendre par 'pied-carré, ce qui est mathématiquement correct mais textuellement problématique. C'est également le choix de Canto-Sperber 1993 dans sa traduction du *Ménon* (82c-84b), signalant toutefois en note le problème.
[11] Bien que cela ne soit pas spécifié, le rectangle ayant pour côtés cet entier et l'unité remplit cette condition. Il s'agit donc dans l'esprit du récit d'un entier, ce qui ne va pas de soi dans l'arithmétique grecque ancienne, l'unité, voire *2*, étant parfois exclus des nombres (selon la définition 2 du livre VII des *Éléments* d'Euclide le nombre est défini comme une 'multitude composée d'unités' ('τὸ ἐκ μονάδων συγκείμενον πλῆθος')). Quoiqu'il en soit, de cette manière on peut toujours associer un rectangle à un entier. Mais inversement, en général, plusieurs rectangles peuvent être associés à un même entier, ainsi à *6*, convient aussi bien le rectangle de côtés *1* et *6* que celui de côtés *2* et *3*.
[12] À savoir *√2*, *√3*, … En effet, on n'est plus ici dans le cadre de la leçon de Théodore qui commence à *3* et s'arrête à *17*. Mais aussi, parce que l'on a quitté la géométrie 'métrique', c'est-à-dire que l'unité de longueur, le pied, utilisé dans la leçon, est abandonnée. Il n'est jamais question d'intermédiaires entre les grandeurs entières et celles irrationnelles i.e. celles correspondant aux rapports d'entiers relativement premiers.

Théodore[13]. Ce qui justifie déjà l'approbation donnée un peu plus loin par Socrate à Théétète[14].

---

[13] Voir aussi la critique de Pappus concernant l'approche métrique de Théodore. Quoique la suite du texte qui se veut explicatif soit assez confuse, peut-être un ajout d'un commentateur ou traducteur ayant mal compris la question, le sens général est clair (Thomson-Junge 1930, §11, p. 74).
[14] Nous ne considérerons pas ici la suite à la leçon de Théodore faite par Théétète et son camarade Socrate, objet d'un prochain article (Ofman 2015).

### 4. Le 'résultat de Socrate-Théétète'

Bien que nous n'étudierons pas ici le travail des deux garçons, il est un point qui nous importe dans le récit qu'en donne Théétète : les côtés des carrés de surfaces des entiers non carrés parfaits sont incommensurables aux côtés des carrés de surfaces entiers carrés parfaits (148b). Elle se traduit, en termes modernes, de la manière suivante :

'Résultat de Socrate-Théétète' :
La racine carrée d'un entier est rationnelle si et seulement si cet entier est un carré parfait[15].

**Remarque.** Il faut souligner, car la confusion a été faite[16], que le 'résultat de Socrate-Théétète' est tout à fait différent du suivant :
La racine carrée d'un entier $n$ est un entier si et seulement si $n$ est un carré parfait.

En effet, ce dernier énoncé qui s'écrit :
$\sqrt{n} = N$ équivaut à : $n = N^2$
est une pure trivialité qui résulte de la définition même de la racine carrée (ou en termes géométriques de la définition de la surface d'un carré en fonction de ses côtés).

Au contraire, le 'résultat de Socrate-Théétète' consiste à montrer l'équivalence :
$\sqrt{n} = p/q$ (où $p$ et $q$ sont des entiers) équivaut à : il existe un entier $r$ tel que $n = r^2$
ou encore, sa contraposée, permettant de préparer une démonstration par l'absurde :
l'entier $n$ n'est pas un carré parfait si et seulement si **quels que soient les entiers $p$ et $q$**, il est impossible que l'on ait : $\sqrt{n} = p/q$ ou encore $nq^2 = p^2$ [17].

La difficulté est l'impossibilité de procéder directement en considérant tous les cas, puisque $p$ et $q$ sont des entiers **quelconques** (donc aussi grands que l'on veut !).

---

[15] Un 'carré parfait' est le carré d'un entier, ainsi $9 = 3^2$ est un carré parfait.

[16] Confusion difficilement évitable eu égard au récit où seuls interviennent les entiers et les irrationnels. De même, dans la solution du problème du doublement du carré dans le *Ménon* (82b-85b), Socrate s'arrange pour n'avoir à considérer que des entiers en évitant les fractions.

[17] En fait dans les deux cas, ce qui importe est l'implication qui va de gauche à droite. Il s'agit dans le premier cas de :
$\sqrt{n}$ est rationnel implique l'existence d'un entier $r$ tel que $n = r^2$ (l'implication inverse est triviale puisque si $n = r^2$, la racine carrée de $n$ est $r$, un entier, donc *a fortiori* elle est rationnelle).
Dans le second cas (la contraposée), cela revient à :
Si l'entier $n$ n'est pas un carré parfait alors il n'existe pas d'entiers $p$ et $q$ tels que l'on ait : $\sqrt{n} = p/q$.
Sa réciproque est évidente puisqu'elle signifie que s'il n'existe pas d'entiers $p$ et $q$ tels que l'on ait :
$\sqrt{n} = p/q$ (donc tels que $n = p^2/q^2$), alors $n$ n'est pas un carré parfait (ce qui suppose outre l'existence de ces entiers $p$ et $q$, que l'on puisse choisir $q = 1$).

## 5. La validité du raisonnement de Socrate-Théétète

Reste à voir la validité du raisonnement des jeunes garçons. Cela est indispensable pour tout énoncé effectué dans un cadre mathématique, et tout aussi crucial pour une enquête devant posséder, dit Socrate, le même degré de certitude qu'une preuve mathématique, et exclure tout ce qui est de l'ordre du vraisemblable ou du probable, au profit du vrai (162e).

Le 'résultat de Socrate-Théétète' est souvent assimilé à la proposition 9 du livre X des *Éléments* d'Euclide, ce qui n'est pas tout à fait exact[18]. C'en est un cas particulier, comme on le notait déjà dans l'Antiquité, ainsi Pappus qui, dans son commentaire du livre X, insistait sur leurs différences (Thomson, Junge 1930, 10-11, p. 73-75)[19].

Il ne manque pas de preuves de ce résultat, toutefois il est un problème délicat : en effet, suivant le récit, Théodore s'est arrêté après avoir étudié, **un à un** ('κατὰ μίαν ἑκάστην'), l'incommensurabilité des côtés de certains carrés de surfaces 3 pieds, 5 pieds. ..., 17 pieds[20]. Il n'a donc certainement pas donné de preuve générale du 'résultat de Socrate-Théétète'.

*A fortiori*, il est tout à fait invraisemblable que deux tout jeunes garçons aient pu le faire, suite à une leçon d'un mathématicien renommé qui, lui, en aurait été incapable. Cette position de bon sens, partagée par les commentateurs modernes[21], ressort du récit platonicien qui ne contient rien approchant d'une quelconque démonstration, ni même de la nécessité d'une démonstration. À la manière dont le conte Théétète, cela apparaît comme une conséquence découlant de source de la leçon, avec la douceur silencieuse d'un flux d'huile qui se répand[22]. Elle paraît en quelque sorte subrepticement introduite par Théétète et son camarade, si ce n'est qu'ils n'en sont nullement conscients.

Le problème est si sérieux que les commentateurs ont beaucoup de mal à le résoudre, et le désaccord est la règle. Certains considèrent qu'il n'y a aucune démonstration chez Théodore, ce qui banalise entièrement cette partie mathématique[23]. L'un des plus longs passages de Platon consacré aux mathématiques serait alors une suite de trivialités, ce qui va à l'encontre aussi bien de la tradition remontant à l'Antiquité que de l'emploi des mathématiques dans ses autres ouvrages[24].
D'autres pensent que pour rendre hommage à Théétète, un collègue et ami, Platon met dans sa bouche une théorie à laquelle il aboutira quelques années plus tard.
C'est oublier la part prise, dans le récit platonicien, par son jeune camarade Socrate, que la tradition n'a retenu d'aucune manière en tant que co-auteur de cette théorie. Il est en outre difficile de comprendre comment l'absence totale de démonstration dans un récit mathématique pourrait être considérée positivement par Platon, alors que, un peu plus loin, Socrate chaudement approuvé par Théétète, insiste sur le caractère essentiel des preuves en mathématiques[25]. Mais les deux objections les plus décisives sans doute, l'une interne et l'autre externe, sont les suivantes.

---

[18] De fait, d'un point de vue strictement mathématique, le résultat du Théétète tel qu'on l'interprète serait plus proche de la proposition VIII.24 des *Éléments* que de la proposition X.9 (cf. aussi supra, note 12).
[19] M. Burnyeat se plaint justement du manque d'intérêt des commentateurs modernes sur cet aspect (Burnyeat 1978, note 60, p. 507).
[20] Cf. supra, note 10.
[21] Ainsi par des auteurs aussi éloignés que Heath, Knorr, Szabo, Caveing ou Burnyeat.
[22] Théodore qualifiait ainsi la facilité d'apprendre de Théétète (144b5).
[23] Ainsi Szabó 1977, p. 63-65
[24] Il suffit de songer à la mathématisation cosmogonique du *Timée* ou à la présentation du nombre 'nuptial' de la *République* (VIII, 546b-c).

Tout d'abord l'insistance répétée de Socrate sur le peu de sérieux qu'il convient d'accorder à ce que peuvent faire des enfants (par exemple, 148c1-4; 165a, d; 166a; 168d-e; 169b9; 200a), voire des adolescents[26], qui ne sont que des jeux comparés à ce que font (ou devraient faire) les adultes. Dans ce cadre, mettre en scène deux jeunes garçons accédant à l'un des plus importants résultats de la géométrie grecque, sans que les mathématiciens de l'époque aient été capables de le démontrer, est inconcevable[27]. D'autre part, il n'est pas raisonnable de penser qu'un lecteur du temps de Platon, contemporain des travaux de Théétète, connaissant les difficultés qu'il a fallu résoudre, aurait pu accepter sans sourciller, dans un ouvrage prétendant à la rigueur de la méthode démonstrative mathématiques, un tel traitement et une telle distorsion de la vérité historique[28]. Moins encore, que son auteur aurait pu le penser, dans un environnement extrêmement polémique, où les écoles philosophiques se concurrençaient férocement[29], et toute maladresse était relevée avec délectation par les adversaires de l'Académie.

---

[25] 'Quant à démonstration et nécessité, il n'y en a d'aucune sorte dans ce que vous dites ; au lieu de cela vous argumentez à coups de 'c'est plausible' : si c'était de cela que voulait se servir Théodore ou quelqu'un d'autre pour parler géométrie, il ne vaudrait seulement rien.' (Narcy 1994) ('ἀπόδειξιν δὲ καὶ ἀνάγκην οὐδ' ἡντινοῦν λέγετε ἀλλὰ τῷ εἰκότι χρῆσθε, ᾧ εἰ ἐθέλοι Θεόδωρος ἢ ἄλλος τις τῶν γεωμετρῶν χρώμενος γεωμετρεῖν, ἄξιος οὐδ' ἑνὸς μόνου ἂν εἴη.' (162e)).

[26] En 168e1-3, Socrate exige de Théodore qu'il accepte d'entrer dans l'analyse de la doctrine de Protagoras afin d'éviter que la discussion puisse être accusée de n'être qu'un 'amusement avec des jeunes gens' ('ὡς παίζοντες πρὸς μειράκια'). Et un peu plus loin, il met en garde à nouveau Théodore sur la possibilité que leur discussion puisse prêter flanc à la critique d'avoir une 'forme qui conviendrait à des enfants' (169b9-c1).

[27] Ce qui conduit certains commentateurs (ainsi Itard 1961, p. 34) à considérer que Platon ne serait qu'un piètre connaisseur des mathématiques de son époque (cf. supra, note 23).

[28] La nécessité pour l'interprétation du texte, de tenir compte des contemporains de Platon à qui il était destiné, est bien souligné par M. Burnyeat (Burnyeat 1978, p. 491).

[29] Songer à l'anecdote rapportée par Diogène Laërce, du poulet déplumé jeté par Diogène le Cynique lors d'une conférence de Platon, pour se moquer de sa définition de l'homme comme bipède sans plumes (Laërce 1965, VI, 2, 40).

## 6. Six critères à vérifier par la démonstration de Théodore

Quant à la leçon de Théodore telle qu'elle est rapportée par Théétète, l'absence d'indications sur la méthode utilisée[30] conduit, comme nous l'avons vu, certains commentateurs à conclure au caractère trivial de cette leçon[31]. Position quelque peu paradoxale dans la mesure où ceux-là mêmes qui l'affirment, utilisent sa supposée insignifiance pour caractériser le récit de Théétète comme une suite de platitudes. Quoiqu'il en soit, ceci est suffisamment improbable (cf. supra, note 21) pour que la plupart admettent l'intérêt de rechercher une méthode de preuve de l'irrationalité des racines carrées de *3, 5, ...* jusqu'à *17*, qui soit en accord avec le récit platonicien[32]

Son omission par Platon (vis-à-vis du lecteur) ou par Théétète (vis-à-vis de Socrate) s'explique naturellement dès que le contexte est pris en compte (cf. note 28, supra). Le problème et sa démonstration étaient bien connus du public mathématiquement cultivé[33], comme celle ayant précédé la théorie générale d'irrationalité. Suivant la tradition, c'est ce même Théétète qui, quelques années plus tard, en sera l'auteur, et c'est elle qui se retrouve, avec quelques modifications, dans les *Éléments* d'Euclide. Aussi, pas plus Théétète pour Socrate, que Platon pour ses contemporains, n'avaient-ils à se soucier de la rappeler. Mais en outre, comme on le verra au paragraphe suivant, la suite des entiers étudiés par Théodore est un indice très fort sur la méthode suivie : elle ne prend sens que si la preuve se fonde sur une certaine propriété des carrés impairs (cf. infra, note 43).

De très nombreuses méthodes de démonstration ont été proposées. Celle qui depuis une cinquantaine d'années a eu le plus de succès parmi les interprètes modernes, se fonde sur ce qu'on appelle, par translittération du grec, 'l'anthyphérèse'. Sans entrer dans les détails, ce qui exigerait un trop long développement, il s'agit d'une généralisation aux grandeurs de l'algorithme arithmétique d'Euclide permettant d'obtenir le Plus Grand Commun Diviseur de deux entiers[34].

---

[30] Du moins en ce qui concerne la construction des carrés de surface *3, 5, ..., 17* pieds, question que se posent nombre de commentateurs, par exemple J. Anderhub, W. Knorr, H. Schmidt, A. Szabo, van der Waerden … *contra* T. Heath qui récuse que le texte implique cette construction (Heath 1981, p. 203, note 2). La position de W. Knorr est ambiguë, car il conclut sur un accord partiel avec la position de Heath (Knorr 1975, p. 74). Mais précisément, il n'est pas question de *cette* construction, Théétète indiquant seulement de manière très vague que Théodore leur a parlé à propos ('περί') de carrés (ou leurs côtés). Quoiqu'il en soit, pour les historiens qui considèrent cette construction comme une donnée d'évidence, ils restent très partagés sur la méthode (cf. par exemple van der Waerden 1963, p. 142-143). Dans la démonstration que nous proposons, cette construction a un caractère secondaire, par contre, elle joue un rôle considérable pour comprendre le type de mathématiques pratiquées par Théodore et l'appréciation qu'elles inspirent à Platon. Si cette question ne nous intéresse pas directement ici, elle sera centrale dans Ofman 2015. Quant à l'argumentation de Heath, elle est essentiellement dirigée contre toute tentative de banaliser la leçon de Théodore.

[31] Par exemple Szabó 1977, p. 66. À l'inverse, M. Burnyeat qui s'oppose radicalement à une telle conclusion, considère que la démonstration est sans importance pour la compréhension du récit platonicien (Burnyeat 1978, p. 505).

[32] Parfois avec un certain scepticisme, ainsi Burnyeat 1978, p. 505.

[33] L'intérêt que le monde grec portait aux questions mathématiques a de quoi surprendre nos contemporains. Ainsi dans une comédie d'Épicharme (vers 500 BCE), on a une remarque sur l'impossible invariance du pair et de l'impair par addition d'une unité (cité Laërce 1965, III, 11). Dans les *Oiseaux* (999-1009), Aristophane, encore un auteur comique, fait une plaisanterie autour de la quadrature du cercle. Ce qui confirme les considérations dans les textes platoniciens, ainsi l'étonnante facilité des interlocuteurs de Socrate à saisir les exemples mathématiques, ou dans les *Lois*, le caractère injurieux que l'Athénien applique à ceux, dont ses interlocuteurs font partie, qui ignorent les questions portant sur l'incommensurabilité, traités de 'pourceaux' (819d-820c), son caractère insultant étant relativisé seulement parce qu'il s'y inclut lui-même. Enfin, on peut encore remarquer l'enthousiasme de jeunes garçons à travailler les mathématiques et de la popularité dont jouit Théodore qui est soulignée au début du *Théétète*.

Cette hypothèse pose cependant de nombreuses difficultés, entre autres par sa complexité et par l'utilisation qu'elle suppose d'un processus infini qui ne se retrouve nulle part dans les ouvrages grecs qui nous sont parvenus. En outre, cette complexité et sa longueur rendent totalement impossible qu'elle puisse être exposée dans le temps d'une leçon de mathématiques (ou même de plusieurs), ce qui est pourtant explicitement affirmé dans le texte platonicien[35].

Au vu de la pléthore de preuves proposées, on pourrait penser que la difficulté principale consiste à déterminer celle utilisée par Théodore. Pourtant, si l'on tient compte des conditions imposées par le récit de Théétète, c'est l'inverse qui est vrai. Il n'en est aucune, en effet, qui soit compatible avec le texte platonicien, lequel impose qu'elle satisfasse six conditions[36] que nous listons ci-dessous[37]. La plupart en vérifient tout au plus une ou deux. Ainsi la démonstration par anthyphérèse respecte la condition ii)[38], mais aucune des autres, pas même la iv), car si d'après les témoignages qui nous sont parvenus elle était sans doute connue et utilisée par les Pythagoriciens, c'était dans un cadre arithmétique et/ou pour approximer certaines grandeurs.

---

[34] Cf. par exemple van der Waerden 1963, p. 145. Il s'agit essentiellement de retirer autant de fois que l'on peut la plus petite grandeur de la plus grande, d'où une troisième grandeur nécessairement plus petite que la deuxième, puis faire de même avec ces deux dernières grandeurs et itérer jusqu'à obtenir, lorsque cela est possible, deux grandeurs dont l'une est le multiple de l'autre, ce qui a lieu si et seulement si les deux grandeurs sont commensurables (*Éléments*, propositions X.2 et X.3).

[35] Nous en ferons une critique plus détaillée dans un prochain travail en cours de rédaction (cf. aussi le cours d'histoire des mathématiques grecques (en anglais), part III, en ligne à l'adresse suivante : http://www.math.jussieu.fr/~ofman). L'argumentation selon laquelle Platon, comme tout auteur, est absolument libre de son texte et des scènes qu'il rapporte ou invente, bien qu'évidemment vraie n'a pas de sens. D'une part cela revient à rendre cette partie purement fictive, donc à lui ôter toute valeur historique, mais surtout c'est encore oublier le contexte historique dans lequel se trouvait l'Ancienne Académie (cf. supra, fin du paragraphe précédent, et en particulier les notes 28 et 29).

[36] Il en est en fait sept, la dernière suit de la toute fin du récit de Théétète qui passe au cas des racines cubiques. Mais elle entre dans le cadre du travail de Théétète et de son camarade Socrate qui, ainsi que nous l'avons dit, n'est pas considéré ici (cf. supra, note 14).

[37] Déjà W. Knorr remarquait la nécessité d'un tel contrôle par le texte, et aboutissait au même constat (cf. Knorr 1975, p. 96-97). Si une preuve vérifie nos conditions, elle vérifie celles de Knorr ; mais l'inverse est faux, en particulier il n'y a rien d'équivalent à notre condition v), pourtant essentielle pour que le récit platonicien soit autre chose qu'un simple conte (ib., p. 193 ; cf. aussi supra note 35).

[38] Ce que ne manquent pas de souligner ses partisans (cf. par exemple van der Waerden 1963, p. 143).

Les six conditions imposées par le texte platonicien sont les suivantes :

i) Débuter à *3* et non pas à *2* ;

ii) Être faite au cas par cas ('κατὰ μίαν ἑκάστην') et non pas donner un résultat général ;

iii) Rendre compte de l'arrêt à *17* ;

iv) Être compatible avec les connaissances mathématiques à l'époque du récit ;

v) Tenir dans le temps d'une leçon destinée à de jeunes garçons ;

vi) Induire à une généralisation abusive à la totalité des entiers.

## 7. Une preuve vérifiant ces six critères

Aussi étrange que cela puisse paraître, on s'est peu interrogé sur la question de savoir quelle était la suite étudiée par Théodore, tant la réponse paraissait évidente.
Pourtant, Théétète rapporte non pas cette suite, mais une abréviation, à savoir *3, 5, ..., 17* autrement dit commençant à *3*, puis *5* et se terminant à *17*, mais rien n'est dit des entiers compris entre *5* et *17*, les historiens des mathématiques débattant de savoir si *17* participe de cette suite[39].

La question que les commentateurs se posent depuis l'Antiquité, est de comprendre pourquoi Théodore débute à *3* plutôt qu'à *2*, la racine carrée de *2* étant irrationnelle, et sa preuve plus simple[40]. La plupart des modernes suivent, quoique parfois avec réticence, l'explication de Hieronymus Zeuthen, suivant laquelle Platon omet ce cas car l'irrationalité de racine carrée de *2* était bien connue[41].

Qu'ils aient raison quant à l'ancienneté de ce résultat est certain. Que ce soit la raison de son absence de la leçon est plus problématique. En effet, une conséquence de la démonstration d'irrationalité de racine de *2* donne beaucoup plus : la question générale de l'irrationalité (des racines carrées d'entiers) est ramenée au cas des entiers impairs, et à prouver qu'ils sont (ou ne sont pas) égaux à un rapport de deux entiers eux-mêmes impairs, auquel cas ils sont rationnels, et dans le cas contraire, ils ne le sont pas[42].

---

[39] W. Knorr note justement que n'importe quelle suite d'entiers différents de *1, 2 et 4,* contenant *3* et *5*, et s'arrêtant à *17*, pourrait convenir.

[40] Dans un article précédent (Ofman 2010), nous en avons donné une démonstration n'utilisant que des résultats très anciens, et qui, contrairement aux preuves usuellement admises, est compatible avec les témoignages textuels. Nous la reprendrons dans ce qui suit. Si elle n'est pas indispensable pour la preuve que nous proposons de la leçon de Théodore, elle la rend plus simple et plus cohérente (cf. infra, note 42). C'est un argument la confortant en tant que démonstration originelle d'incommensurabilité. Du point de vue de la démonstration par *anthyphérèse*, le cas de *2* est le plus simple, ce qui ajoute un argument à l'encontre de l'utilisation de cette méthode par Théodore. En effet :
- Considérer que l'irrationalité de racine carrée de *2* était également prouvée par anthyphérèse conduit à revoir les datations concernant son apparition dans les mathématiques grecques, et donc à de nouvelles difficultés voire des pétitions de principe. Théodore étant contemporain de ce résultat, s'il n'en est pas l'auteur, son apport est considérablement réduit et aussi bien l'intérêt de sa leçon. En outre, qu'il en soit ou non l'auteur, l'absence du cas *2* est difficilement explicable. Et surtout, ceci va directement à l'encontre des témoignages textuels que nous avons sur ces questions.
- Si au contraire on suppose que le cas de *2* était connu et prouvé de manière différente, Théodore étant alors le premier à l'utiliser pour prouver d'autres cas d'irrationalité, cela ne rend pas mieux compte de son absence de la leçon. En effet, la construction *anthyphérétique* étant la plus simple pour *2*, un enseignant ne pourrait que débuter par lui, pour ensuite seulement passer à une généralisation plus complexe. Dans tous les cas, il aurait pour le moins mentionné cette possibilité de validité universelle de sa méthode, l'absence d'une telle indication mettant d'ailleurs en cause la crédibilité du récit platonicien tout entier (cf. supra, note 28).

[41] Considérée comme un pis-aller, faute d'une autre explication, ainsi Burnyeat 1978, p. 502-503.

[42] Pour plus de détails, cf. Ofman 2010, p. 118.

Autrement dit, pour tout entier *n* impair, on est conduit à rechercher s'il existe (ou pas) des impairs *p* et *q* tels que : $\sqrt{n} = p/q$ [43].

La suite donnée par Théétète : *3, 5, ..., 17* est donc simplement une abréviation de la suite naturelle des impairs compris entre *3* et *17*. Pour une telle suite, toute explication supplémentaire était superflue, que ce soit pour Socrate ou pour le lecteur contemporain de Platon.

Reste à montrer la rationalité ou l'irrationalité des racines carrées de ces entiers.
La preuve que nous proposons est fondée sur un résultat déjà connu des anciens Pythagoriciens concernant les carrés des impairs[44], et étroitement liée à la suite étudiée par Théodore[45]. Il s'agit du

**Résultat du reste**' : La division par *8* d'un impair au carré a pour reste *1*.

En particulier, s'il est différent de l'unité, il doit être plus grand que *8*.
En termes modernes, on dirait que tout carré (parfait) impair est égal à *1 modulo 8*.

Ce résultat est conséquence de la

**Remarque très simple**. Si on ajoute *8* à un entier quelconque *n*, le reste sa division par *8* ne change pas. Il en est donc de même si l'on ajoute un multiple de *8* à *n* (i.e. pour tous entiers *m* et *n*, les entiers *n* et *n+8m* divisés par *8* ont même reste).

Sa démonstration est extrêmement simple. En symbolique moderne, il suffit d'écrire :
*(2n+1)$^2$ = 4n$^2$ + 4n + 1 = 4n(n+1) + 1*.
Le produit *n(n+1)* est pair (puisque soit *n*, soit *n+1* l'est) et donc  *4n(n+1)* est un multiple de *8*. Aussi (cf. la remarque très simple ci-dessus) *4n(n+1) + 1* divisé par *8* a pour reste *1*.

---

[43] Si cela peut se déduire plus ou moins facilement de la plupart des démonstrations, c'est une conséquence pratiquement évidente de celle donnée dans Ofman 2010. Que l'on soit ramené au cas des impairs ne signifie pas que l'on sache résoudre la question de l'irrationalité des racines carrées des entiers pairs. Ainsi, comme conséquence de la démonstration, on obtient immédiatement l'irrationalité de la racine carrée de *6*. Par contre celle de *12 = 4 × 3* est égale au double de celle de *3* (car $\sqrt{12} = \sqrt{4 \times 3} = 2\sqrt{3}$) ; elle est donc rationnelle (ou irrationnelle) si et seulement si la racine carrée de *3* l'est.
[44] Par exemple Plutarque, *Questions Platoniques*, II, 24, 1003f. C'est Jean Itard (Itard 1961, p. 34-36) qui a attiré l'attention sur cette proposition, reprise par divers historiens des mathématiques, dont M. Caveing, W. Knorr, J. Vuillemin. De fait, la seule contemplation des tableaux des carrés impairs conduisait à cette déduction, sa démonstration étant d'une grande simplicité même dans le cadre d'une arithmétique très primitive, ce qui conduit à la possibilité d'un résultat déjà connu des calculateurs babyloniens ou égyptiens (ib., p. 34).
[45] Il faut souligner une fois encore que la difficulté qui se pose à nous concernant la démonstration employée par Théodore pour obtenir les résultats exposés dans le récit, n'en était pas une pour les lecteurs de Platon, la suite étudiée signant d'ailleurs la méthode utilisée.

La démonstration géométrique est tout aussi simple, et pouvait prendre la forme suivante[46] :

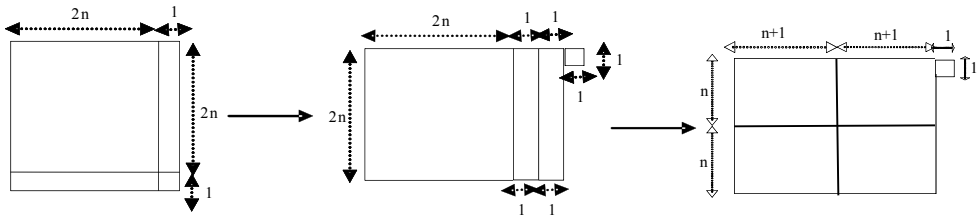

Il reste à voir comment ce résultat, qui porte uniquement sur les entiers et non sur (ce que nous appelons) les rationnels, peut être utilisé (cf. supra, remarque du §4).

D'après ce que nous avons noté au début de ce paragraphe (cf. supra, note 42), pour que la racine carrée d'un entier $n$ impair soit rationnel, il faut (et il suffit) qu'il existe des entiers **impairs** $p$ et $q$ en sorte que l'on ait :
$\sqrt{n} = p/q$.
Si de tels entiers existent alors la racine carrée de $n$ est rationnelle, sinon elle est irrationnelle.

En élevant au carré les deux termes de cette égalité, on obtient :
$n = p^2/q^2$
ou encore
$nq^2 = p^2$.

On applique alors le 'résultat du reste' tout d'abord à $q^2$ (carré d'un impair), et on obtient :
$nq^2 = n(8k+1) = 8nk + n$.

Puisque $nq^2 = p^2$, c'est le carré d'un impair et on peut appliquer le 'résultat du reste' cette fois à $nq^2 = 8nk + n$
et le reste de la division par $8$ de $8nk + n$ est donc égal à $1$.

Mais (cf. la remarque très simple ci-dessus) le reste de la division par $8$ de :
$n(8k+1) = 8nk + n$, est égal au reste de $n$ divisé par $8$, et l'on obtient finalement :

le reste de $n$ divisé par $8$ est égal à $1$           (\*).

**Remarque**. La démonstration algébrique ci-dessus est évidemment anachronique. Nous en donnons ci-dessous une démonstration géométrique où le résultat se *voit* sur les figures.

---

[46] Sur l'attribution aux Pythagoriciens anciens de ces résultats, cf. par exemple Proclus (Proclus 1992), en particulier son commentaire sur les premières propositions du livre I des *Éléments* d'Euclide. La relation $(m+n)^2 = m^2+2mn+n^2$ ($m$ et $n$ entiers donnés), dont la connaissance est nécessaire pour les plus anciennes démonstrations connues de cas particuliers du théorème de Pythagore, est déjà utilisée dans des tablettes mésopotamiennes du 2[ème] millénaire BCE, ainsi que le montrent par exemple celles de la Collection Schøyen (cf. Friberg 2007, p. 42-51).

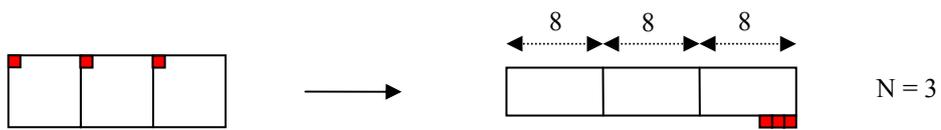

N = 3

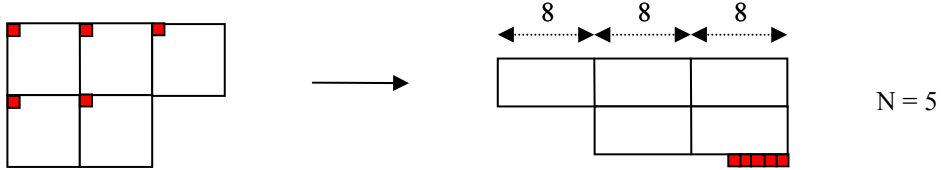

N = 5

..................................................................................

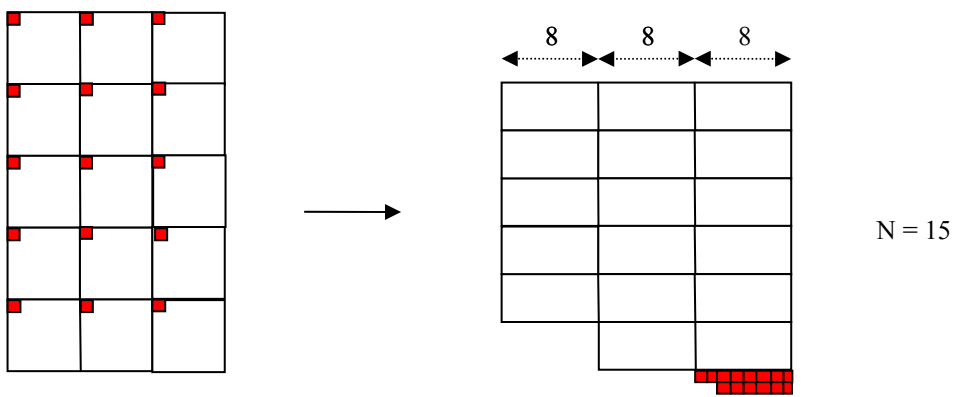

N = 15

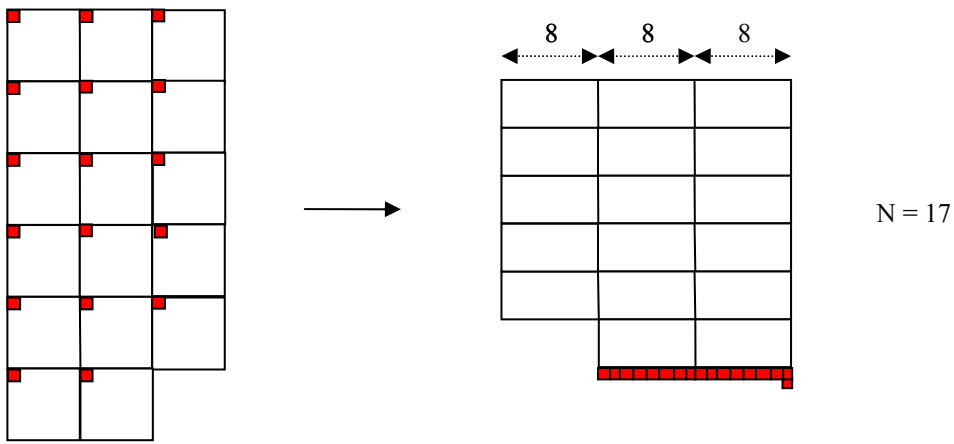

N = 17

Dans la figure, les parties colorées représentent une unité (en langage de Théodore un carré d'un pied), et les parties blanches des (surfaces) multiples de *8*. La flèche permettant de passer de la figure de gauche à celle de droite est le réarrangement des carrés tronqués (i.e. des carrés dont on a ôté une unité, donc des 'gnomons'[47]) en des rectangles dont un côté est *8* unités. Il

---

[47] Instrument en forme d'équerre utilisé en Grèce ancienne dans toutes sortes de situations, en particulier pour des calculs astronomiques.

*apparaît* alors immédiatement que le reste de la division par *8* de *3* fois n'importe quel carré impair a pour reste *3* (et non pas *1*), etc…

Considérons alors le tableau où l'on a reporté les entiers de la suite de Théodore et leur reste dans la division par *8* :

| *n* | reste de la division de *n* par *8* |
|---|---|
| 3 | 3 |
| 5 | 5 |
| 7 | 7 |
| 9 | 1 |
| 11 | 3 |
| 13 | 5 |
| 15 | 7 |
| 17 | 1 |

Les trois premiers nombres (*3, 5, 7*) ne vérifient pas (*), aussi d'après le 'résultat du reste', leurs racines carrées ($\sqrt{3}$, $\sqrt{5}$, $\sqrt{7}$) ne sont pas rationnelles.
Le premier problème se rencontre avec *9*, mais *9* est précisément un carré parfait (*9 = 3²*), sa racine carrée (à savoir *3*) est entière, donc rationnelle.

De même, les trois entiers suivants (*11, 13, 15*) ne vérifient pas (*), et leurs racines carrées ($\sqrt{11}$, $\sqrt{13}$, $\sqrt{15}$) sont encore irrationnelles.

Reste alors *17* qui divisé par *8* a pour reste *1*. Mais cette fois *17* n'est pas un carré parfait, et les deux garçons d'en conclure naturellement que **sa racine carrée est irrationnelle**.
Pourtant, si on peut certainement conclure que $\sqrt{17}$ **n'est pas entier**[48], on ne peut rien dire **quant à sa rationalité ou pas**.

Ainsi Théodore s'est arrêté là sans donner de conclusion. Soit afin de laisser les élèves trouver le résultat par eux-mêmes (thèse retenue par A. Szabó[49]), soit plus probablement afin qu'ils reviennent suivre la séance suivante[50].

Cette manière de procéder n'est toutefois pas sans danger Ainsi les deux jeunes garçons ont certes remarqué que la question pouvait s'étendre à une infinité d'entiers[51] (147d), et que la démonstration au cas par cas de Théodore ne pouvait y répondre, d'où leur tentative d'obtenir une théorie générale. C'est une nouvelle justification[52] du compliment qui leur est adressé un peu plus loin par Socrate (147e8, 148b7). Toutefois, la démonstration de Théodore[53] ne mettant pas en évidence qu'il y ait là un problème majeur, l'arrêt paraît incompréhensible : Théodore, dit Théétète, 's'était, je ne sais pourquoi, arrêté là [à 17]' (147d7)[54]. Et les deux jeunes garçons d'en conclure à une alternative identique aux cas précédemment étudiés :
- ou bien on peut appliquer le résultat de la division par *8*,
- ou bien lorsqu'on ne le peut pas, nouvelle alternative :
    - ou bien le nombre est un carré parfait (cas de *9*) et sa racine carrée est rationnelle
    - ou bien le nombre ne l'est pas (cas de *17*) et sa racine carrée est irrationnelle.

---

[48] En effet, *17* est plus grand que *16 = 4²* et plus petit que *25 = 5²*, or il n'existe pas d'entier entre *4* et *5*.
[49] Szabó 1977, p. 92.
[50] Les enseignants étaient payés en fonction de leur audience, et lorsque Socrate est confronté à des enseignants (souvent des 'sophistes'), les questions d'argent ne sont jamais très loin. Ainsi au tout début de la discussion, la première chose dont parlent Théodore et Socrate porte sur l'argent, plus précisément celui de Théétète (140c6-d3). Et plus loin, lorsqu'il s'agit de discuter des thèses de Protagoras, l'ami de Théodore, Socrate rappelle que celui-ci exigeait de ses élèves 'd'énormes honoraires' ('μεγάλων μισθῶν', 161d).
[51] Cela est déjà évident dès que l'on sait que $\sqrt{2}$ est irrationnel, puisque pour tout entier *n*, on a également :
$\sqrt{(2n^2)} = n\sqrt{2}$ est irrationnel. Cette égalité utilisée par Socrate dans le *Ménon* (83a-c), pour montrer au jeune esclave comment doubler la surface d'un carré, était donc très anciennement connue. Elle apparaît d'ailleurs évidente dès que l'on trace la figure.

Et si Théodore examine bien le cas *17*, c'est pour le laisser ouvert.

Cette démonstration fondée sur l'utilisation du 'résultat du reste', vérifie bien les six critères que nous avons listés. En effet :
- Elle explique le début à *3* et l'arrêt à *17*.
- Elle est faite au cas par cas et n'utilise que des connaissances arithmétiques très anciennes, sans doute adaptées par Théodore au cadre de l'irrationalité.
- Elle tient aisément dans le délai d'une heure.
- Elle explique enfin la précipitation des jeunes garçons à conclure sur un énoncé qu'ils ne sont aucunement en état de prouver (cf. supra, fin du §5).

---

[52] Cf. supra, note 12.
[53] Démonstration qui utilise exclusivement des grandeurs entières et non des grandeurs commensurables plus générales, ce qui est précisément relevé par Pappus (§10, p. 73 ; cf. aussi supra, note 18).
[54] Il s'agit de la traduction de 'ἐν δὲ ταύτῃ πως ἐνέσχετο' qui est très discutée (cf. Burnyeat 1978, p. 513-514 contra Knorr 1975, p. 62). Ainsi dans Narcy 1994, cela est rendu par 'quelque chose l'a arrêté [là]', pour d'autres encore il s'agirait d'un arrêt purement aléatoire ce qui n'est pas vraiment conforme au texte grec.

## 8. Les conclusions

### i) Le caractère aporétique de la démonstration de Théodore

On considère généralement que la 'partie mathématique du Théétète' rend compte des découvertes de Théodore sur l'irrationalité, et prophétiquement de la future théorie des irrationnels, telle qu'on la trouve au livre X des *Éléments* d'Euclide traditionnellement associé à Théétète devenu mathématicien. Ce texte est alors très difficile à relier au reste de l'ouvrage platonicien.

La tradition exégétique souligne en effet que celui-ci est essentiellement aporétique, et aboutit à une série d'échecs. Toutes les conclusions sont, au grand dam de Théodore, de Théétète[55] et sans doute aussi bien du lecteur, minutieusement réfutées au fur et à mesure par Socrate. Tant est qu'il est bien difficile de préciser la position de Platon (ou de Socrate) vis-à-vis des définitions successives qui sont proposées.

Par contre, le récit mathématique ouvrirait sur l'un des plus grands succès des mathématiques grecques classiques, la théorie des irrationnels. Dès lors, cette partie, loin d'éclairer le dialogue, serait facteur de confusion[56]. Suivant l'analyse que nous avons proposée, ce début de la 'partie mathématique du *Théétète*' s'inscrit, au contraire, de manière cohérente dans l'ouvrage.

Certes, les résultats de Théodore permettent d'obtenir un très grand nombre de réponses. En langage moderne, on dirait dans 80% des cas[57]. Du point de vue 'effectif' (au sens de l'algorithmique moderne), le processus est performant, et Théodore (et Théétète) a toutes les raisons d'être fier de lui. Toutefois, du point de vue mathématique, aussi généreux qu'il soit, cela ne permet pas de caractériser la rationalité/irrationalité des racines carrées des entiers comme le voudraient les deux garçons (147e1) suivant la demande de Socrate (147c1, 146d3-4, 146e7-9, 147e1, 148d4-7). La mathématique grecque ultérieure ne reprendra d'aucune manière le 'résultat du reste' dans ce cadre. Il sera abandonné au profit de la théorie des proportions et des entiers relativement premiers.

On est bien dans une aporie, dont on sort en changeant de méthode, tout comme indubitablement, il s'agira de le faire pour répondre à la question posée par Socrate : 'qu'est-ce que la science ?'.

Mais inversement, suivant cette analogie, l'enquête sur la science n'est pas simplement une aporie, et la dernière définition indiquée par Théétète, de la science comme 'opinion vraie accompagné d'un raisonnement/d'une définition ('λόγος')' serait plus pertinente pour définir la science que la liste donnée par Théétète au tout début de l'entretien (cf. supra, §2). Poursuivant ce parallélisme, la dernière définition ne serait certes pas celle de 'la science', mais du moins de la plupart 'des sciences'. Et le point essentiel sur lequel achopperait l'enquête serait alors l'oubli d'une division nécessaire préalable entre 'les sciences'[58]. La méthode définitionnelle de divisons sera rétablie par l'Étranger dans le *Sophiste*, présenté comme la suite du *Théétète*, et le *Politique*[59].

---

[55] Entre autres, 157c4-6, 161a4-9, 165d1.
[56] Ainsi B. van der Waerden, quoique pas toujours très rigoureux lorsqu'il présente le texte platonicien (par ex. van der Waerden 1963, p. 141), note pourtant qu'elle apparaît comme 'une pièce rapportée (…) qui n'y a pas vraiment sa place' (ib., p. 166 ; et aussi p. 142).
[57] Ou encore que la probabilité d'obtenir une réponse est de 80%. Au sens où, étant donné un nombre pris au hasard, on a plus 8 chances sur 10 pour savoir si sa racine carrée est irrationnelle. Ainsi dans la suite de *1* à *17* (compris), il n'y a qu'un seul cas auquel l'on ne puisse répondre.

**ii) La critique de Théodore par Platon**

Cette critique a lieu sur deux plans.
 D'une part, sur la possibilité de faire des mathématiques à la manière de Théodore, non par des raisonnements, mais, comme le rapporte Théétète, par des dessins ('ἔγράφε') (147d2), et au cas par cas ('κατὰ μίαν ἑκάστην', 147d6). Or dès que l'infini apparaît ou qu'une démonstration par l'impossible est nécessaire, cette procédure est vouée à l'échec. Selon le récit, la leçon de Théodore apparaît même en-deçà de ce que peuvent faire les deux jeunes garçons qui passent aussitôt à un raisonnement sur les entiers en général, abandonnant les mesures physiques données en 'pieds' dans la leçon. Ses mathématiques se trouvent ainsi reléguées en quelque sorte au rang d'un jeu infantile.

Platon n'est pas plus tendre concernant sa manière d'enseigner. Le dessin dissimule les difficultés, tout comme un vêtement dissimule les corps (165a1) : ce qui ne peut être dessiné, et donc apparaître dans des figures, est simplement ignoré.
D'un point de vue mathématique, cela se traduit de la manière suivante. Trois cas sont à distinguer puisqu'une racine carrée peut être :
- soit entière
- soit non entière mais rationnelle
- soit enfin irrationnelle.

Et deux seulement (les premier et troisième cas) sont considérés par les jeunes garçons.
Ce sont, il est vrai, les seuls effectivement possibles, les seuls donc que l'on puisse faire apparaître et *a fortiori* dessiner. Mais la démarche démonstrative mathématique ne peut s'en contenter. Il faut le prouver, et seul un raisonnement peut le faire. Ce sont précisément les dessins qui induisent en erreur les jeunes garçons, amenés à croire évident un résultat qui va changer la mathématique grecque[60].

Ce jugement par Platon sur les mathématiques à la Théodore, et par ricochet celles du jeune Théétète, du moins au temps où il était son élève, implique un changement de perspective sur l'ouvrage. Contrairement à l'interprétation usuelle, loin de se livrer à une louange de Théodore qui aurait été son maître, et rendre hommage à Théétète en tant que collègue et ami, Platon condamne les mathématiques qu'ils (re)présentent ici.

---

[58] Socrate avait déjà imposé une division entre les 'sciences' et les 'techniques artisanales' que Théétète avait regroupées ensemble (cf. supra, note 5).
[59] Dans un autre ouvrage, Socrate se révèle 'amoureux des divisions ('διαιρέσεων')' (*Phèdre*, 266b4).
[60] La théorie des proportions, c'est-à-dire essentiellement le livre V des *Éléments* d'Euclide, attribuée traditionnellement à Eudoxe, ne prend vraiment son sens que lorsqu'on sort du cadre de la commensurabilité des grandeurs, en termes modernes des rationnels.

Il n'est pas sans intérêt de remarquer que l'étude mathématique de cette partie[61] rejoint le point de vue développé par Michel Narcy (Narcy 1994, introduction, particulièrement pp. 40-69). L'auteur procède à partir d'une analyse purement textuelle et philosophique, par comparaison à d'autres textes de Platon sur les mathématiques, en particulier ceux de la *République*. C'est en quelque sorte une confirmation extérieure de notre analyse fondée sur les questions mathématiques posées par ce passage.

---

[61] Comme le lecteur l'aura sans doute remarqué, nous avons évité ici toutes les questions philologiques ou de traduction. Nous reviendrons sur elles dans Ofman 2015.

# Bibliographie

BURKERT, Walter.1972. *Lore and Science in Ancient Pythagoreanis.* Translated by Edwin Minar. Harvard University Press (1962).

BURNYEAT, Myles. 1978. « The philosophical sense of *Theaetetus*' Mathematics ». *ISIS*, 69 (n° 249), 489-513.

CAVEING, Maurice. 1997. *La constitution du type mathématique de l'idéalité dans la pensée grecque*. 3 volumes. Presse univ. Septentrion.

CHAPPELL, Timothy. 2004. *Reading Plato's Theaetetus*. Translation, notes and commentary, Academia Verlag.

CANTO-SPERBER, Monique. 1993. Platon Ménon. Traduction, introduction et notes. Flammarion.

DIÈS, Auguste. 1976. Platon *Théétète*. Les Belles Lettres (1926).

FOWLER, Harold. 1996. Plato *Theaetetus*. Harvard Univ. Press (1921).

FRIBERG, Jöran. 2007. *A Remarkable Collection of Babylonian Mathematical Texts: Manuscripts in the Schøyen Collection Cuneiform Texts* I, Springer.

HEATH, Thomas. 1956. *The thirteen books of Euclid's Elements*. 3 volumes. Dover (1908).

HEATH, Thomas. 1981. *A History of Greek Mathematics*. Dover (1921).

ITARD, Jean. 1961. *Les livres arithmétiques d'Euclide*. Hermann.

KNORR, Wilbur. 1975. *The Evolution of the Euclidean Elements*. Reidel.

LAËRCE, Diogène. 1965. *Vies, doctrines et sentences des philosophes illustres*. Garnier Flammarion.

McDOWELL, John. 1973. Plato *Theaetetus*. Translation, notes and commentary. Clarendon Press.

NARCY, Michel. 1994. Platon *Le Théétète*. Traduction, introduction et notes. Flammarion.

OFMAN, Salomon. 2010. « Une nouvelle démonstration de l'irrationalité de racine carrée de *2* d'après les *Analytiques* d'Aristote », *Philosophie antique,* n° 10, 81-138.

OFMAN, Salomon. 2015. « Un *logos alogos* : le « théorème de Théétète » (*Théétète*, 147d8-148b2) ». *Philosophie antique*, à paraître.

HERMANN Schmidt, 1877. *Kritischer Kommentar zu Plato*. Teubner.

PROCLUS. 1992. *A commentary on the first book of Euclid's* Elements. Translated by Glenn Morrow. Princeton Univ. Press.

SZABÓ, Àrpád. 1977. *Les débuts des mathématiques grecques*. Traduction par M. Federspiel. Vrin (1969).

THOMSON, William, JUNGE, Gustav. 1930. *The Commentary of Pappus on Book X of Euclid's Elements*. Translation and notes. Harvard University Press.

VITRAC, 1990-2001. Euclide *les Éléments*. Traduction, notes et commentaires. PUF, 1990-2001.

VUILLEMIN, Jules. 2001. *Mathématiques pythagoriciennes et platoniciennes*. Blanchard.

Van der WAERDEN, Bartel. 1963. *Science Awakening*. Translation Arnold Dresden. Science Editions (1950).